\documentclass{svjour3}                     
\smartqed  
\usepackage{graphicx}

\usepackage{amssymb}
\usepackage{amsfonts}
\usepackage{latexsym}
\usepackage{amscd}
\usepackage[mathscr]{euscript}
\usepackage{xy} \xyoption{all}

\usepackage{amsmath}

\newcommand\limind{\mathop{\oalign{\hfill $\rm lim$\hfill\cr$\longrightarrow$\cr}}}

\newcommand{\N}{{\mathbb{N}}}
\newcommand{\Z}{{\mathbb{Z}}}

\newcommand{\uloopr}[1]{\ar@'{@+{[0,0]+(-4,5)}@+{[0,0]+(0,10)}@+{[0,0] +(4,5)}}^{#1}}
\newcommand{\uloopd}[1]{\ar@'{@+{[0,0]+(5,4)}@+{[0,0]+(10,0)}@+{[0,0]+ (5,-4)}}^{#1}}
\newcommand{\dloopr}[1]{\ar@'{@+{[0,0]+(-4,-5)}@+{[0,0]+(0,-10)}@+{[0, 0]+(4,-5)}}_{#1}}
\newcommand{\dloopd}[1]{\ar@'{@+{[0,0]+(-5,4)}@+{[0,0]+(-10,0)}@+{[0,0 ]+(-5,-4)}}_{#1}}

\newcommand{\luloop}[1]{\ar@'{@+{[0,0]+(-8,2)}@+{[0,0]+(-10,10)}@+{[0, 0]+(2,2)}}^{#1}}

 \journalname{Algebras and Representation Theory}
\begin{document}

\title{Regularity conditions for arbitrary Leavitt path algebras
}



\author{Gene Abrams        \and
        Kulumani M. Rangaswamy 
}


\institute{Gene Abrams and Kulumani M. Rangaswamy \at
              Department of Mathematics \\
              University of Colorado at Colorado Springs \\
              Colorado Springs, Colorado, 80933  U.S.A. \\
              \email{abrams@math.uccs.edu  \ \   krangasw@uccs.edu} \\
              Corresponding author:  G. Abrams  \  7192623182(v)  \  7192623605(f)}


\date{Received: date / Accepted: date}

\maketitle

\begin{abstract}
We show that if $E$ is an arbitrary acyclic graph then the Leavitt path algebra $L_K(E)$ is locally $K$-matricial;
that is, $L_K(E)$ is the direct union of subalgebras, each isomorphic to a finite direct sum of finite matrix rings over the field $K$.  (Here an {\it arbitrary} graph means that neither  cardinality conditions nor graph-theoretic conditions (e.g. row-finiteness) are imposed on $E$. These unrestrictive conditions are in contrast to the hypotheses used in much of the literature on this subject.) As a consequence we get our main result, in which we show that the following conditions are equivalent for an arbitrary graph $E$:  (1)  $L_K(E)$ is von Neumann regular.  (2)  $L_K(E)$ is $\pi$-regular.  (3)  $E$ is acyclic.   (4)  $L_K(E)$ is locally $K$-matricial.   (5)  $L_K(E)$ is strongly $\pi$-regular.    We conclude by showing how additional regularity conditions (unit regularity,  strongly clean) can be appended to  this list of equivalent conditions.

\end{abstract}

\bigskip
\bigskip

In this article we investigate Leavitt path algebras $L_K(E)$ over
an arbitrary  directed graph $E$ and an arbitrary field $K$. Here we require
no restriction on either the size of the graph (i.e., vertex and edge sets may be
of any cardinality), or on graph-theoretic constraints (i.e.,  no row-finiteness conditions
are assumed).  Our goal is to classify in ring-theoretic terms the Leavitt path algebras
of the form $L_K(E)$ where $E$ is acyclic.   A useful tool in our study is a
construction presented in Proposition \ref{KeyLemma}, which enables us to realize $L_K(E)$
as a  directed union of subalgebras, each of which is a Leavitt path algebra of a suitable finite
graph. An effective use of Proposition \ref{KeyLemma} leads to our main result, Theorem \ref{VNRegular},
in which we establish the equivalence of these conditions for an arbitrary directed graph
$E$ and field $K$:
(1) $L_K(E)$ is  von Neumann regular.
(2) $L_K(E)$ is $\pi$-regular.
(3) $E$ is acyclic.
(4) $L_K(E)$ is locally $K$-matricial (i.e., $L_K(E)$
is a direct union of subalgebras, each of which is isomorphic to
a finite direct sum of finite matrix rings over $K$).
(5)  $L_K(E)$ is strongly $\pi$-regular.
We conclude the article by discussing various additional ring-theoretic conditions which in the context of Leavitt path algebras are equivalent to $E$ being acyclic.

We begin by giving a terse reminder of the germane definitions.  For a more complete description and discussion, see e.g. \cite{AA3} or \cite{G}.    A \emph{(directed) graph} $E=(E^0,E^1,r,s)$ consists of two  sets $E^0,E^1$ and functions
$r,s:E^1 \to E^0$.  (The sets $E^0$ and $E^1$ are allowed to be of arbitrary cardinality.)  The elements of $E^0$ are called \emph{vertices} and the elements of $E^1$ \emph{edges}.  A \emph{path} $\mu$ in a graph $E$ is a sequence of edges
$\mu=e_1\dots e_n$ such that $r(e_i)=s(e_{i+1})$ for $i=1,\dots,n-1$. In this case, $s(\mu):=s(e_1)$ is the
\emph{source} of $\mu$, $r(\mu):=r(e_n)$ is the \emph{range} of $\mu$, and $n$ is the \emph{length} of $\mu$.
 We view the elements of $E^0$ as paths of length $0$.  If $\mu = e_1...e_n$ is a path in $E$, and if $v=s(\mu)=r(\mu)$ and $s(e_i)\neq s(e_j)$ for every $i\neq j$, then $\mu$ is called a
\emph{cycle based at} $v$.
If
$s^{-1}(v)$ is a finite set for every $v\in E^0$, then the graph $E$ is called \emph{row-finite}.

\medskip

\begin{definition}\label{definition}  {\rm Let $E$ be any directed graph, and $K$ any field. The
{\em Leavitt path $K$-algebra} $L_K(E)$ {\em of $E$ with coefficients in $K$} is  the $K$-algebra generated by a set
$\{v\mid v\in E^0\}$ of pairwise orthogonal idempotents, together with a set of variables $\{e,e^*\mid e\in E^1\}$,
which satisfy the following relations:

(1) $s(e)e=er(e)=e$ for all $e\in E^1$.

(2) $r(e)e^*=e^*s(e)=e^*$ for all $e\in E^1$.

(3) (CK1) $e^*e'=\delta _{e,e'}r(e)$ for all $e,e'\in E^1$.

(4) (CK2)  $v=\sum _{\{ e\in E^1\mid s(e)=v \}}ee^*$ for every vertex $v\in E^0$ having $1\leq |s^{-1}(v)| < \infty$.
}

\end{definition}

\medskip


 For any
$F\subseteq E^1$ the set
$\{e^*\mid e\in F\}$ will be denoted by $F^*$. We let $r(e^*)$ denote $s(e)$, and we let $s(e^*)$ denote $r(e)$.
If $\mu = e_1 \dots e_n$ is a path, then we denote by $\mu^*$ the element $e_n^* \dots e_1^*$ of $L_K(E)$.


Many well-known algebras arise as the Leavitt path algebra
of a graph.  For instance, the classical Leavitt
algebras $L_K(1,n)$ for $n\ge 2$ arise as the algebras $L_K(R_n)$ where $R_n$ is the ``rose with $n$ petals" graph
 described in Example \ref{EsubFexample} below.  (See e.g. \cite[Section 3]{AALP}.)  Also, for each $n\in  \N = \{1,2,...\}$, the full matrix ring ${\rm M}_n(K)$  arises as the Leavitt path
algebra of the oriented $n$-line graph
$$\xymatrix{{\bullet}^{v_1} \ar [r] ^{e_1} & {\bullet}^{v_2} \ar [r] ^{e_2} & {\bullet}^{v_3} \ar@{.}[r] &
{\bullet}^{v_{n-1}} \ar [r] ^{e_{n-1}} & {\bullet}^{v_n}} $$ while the Laurent polynomial ring $K[x,x^{-1}]$ arises as
the Leavitt path algebra of the ``one vertex, one loop" graph
$$\xymatrix{{\bullet}^{v} \ar@(ur,dr) ^x}$$

\medskip

A (possibly nonunital) ring $R$ is called a {\it ring with local units} in case for each finite subset $S\subseteq R$ there is an idempotent $e\in R$ with $S\subseteq eRe$. If $E$ is a graph for which $E^0$ is finite then we have $\sum _{v\in E^0} v$ is the multiplicative identity in $L_K(E)$; otherwise, $L_K(E)$ is a ring with a set of local units
consisting of sums of distinct vertices. Conversely, if $L_K(E)$ is unital, then $E^0$ is finite. $L_K(E)$ is a
${\mathbb Z}$-graded $K$-algebra, spanned as a $K$-vector space by $ \{pq^* \mid p,q$ are paths in $E\}$.  (Recall that the elements of $E^0$ are viewed as paths of length $0$, so that this set includes elements of the form $v$ with $v\in E^0$.)   In particular,
for each $n\in \mathbb{Z}$, the degree $n$ component $L_K(E)_n$ is spanned by elements of the form  $\{pq^* \mid {\rm
length}(p)-{\rm length}(q)=n\}$.  The degree of an element $x$, denoted $deg(x)$, is the lowest number $n$ for which
$x\in \bigoplus_{m\leq n} L_K(E)_m$. The $K$-linear
extension of the assignment $pq^* \mapsto qp^*$ (for $p,q$ paths in $E$) yields an involution on $L_K(E)$, which we
denote simply as ${}^*$.

A subgraph $G$ of a graph $E$ is called {\it complete} in case, for each $v\in G^0$ having $1\leq |s_G^{-1}(v)| < \infty$, we have  $s^{-1}_G(v) = s^{-1}_E(v)$.  (In other words, a subgraph $G$ of $E$ is complete if, whenever $v\in G^0$ emits a nonzero, finite number of edges in $G$, then necessarily the subgraph $G$ contains all of the edges in $E$ emitted by $v$.)   The natural inclusion map $L_K(G) \mapsto L_K(E)$ is a ring homomorphism precisely when $G$ is a complete subgraph of $E$, so that complete subgraphs of $E$ naturally give rise  to subalgebras of $L_K(E)$.     One of our main objectives in this article is to show how to construct subalgebras of $L_K(E)$ which need not arise in this way.  This in turn will allow us to describe algebras of the form $L_K(E)$ as unions of subalgebras possessing various ring-theoretic properties, even in situations where $E$ lacks complete subgraphs possessing corresponding graph-theoretic properties.  We achieve this objective in Proposition \ref{KeyLemma}.  The construction  is based on an idea  presented by Raeburn and Szyma\'{n}ski in   \cite[Definition 1.1]{RS}; this work was brought to our attention by E. Pardo.

\medskip

\begin{definition} \label{EsubFdefinition} {\rm
Let $E$ be a graph, and let $F$ be a finite set of edges in $E$.   We define $s(F)$ (resp. $r(F)$) to be the sets of those vertices in $E$ which appear as the source (resp. range) vertex of at least one element of $F$.  We define a graph $
E_{F}$ as follows:
$$E_{F}^{0}=F \ \cup \ (r(F)\cap s(F)\cap s(E^{1}\backslash F)) \ \cup \ (r(F)\backslash s(F)),$$
$$E_{F}^{1}=\{(e,f)\in F\times E_{F}^{0}\mid r(e)=s(f)\} \ \cup \  [\{(e,r(e))\mid e\in F \mbox{ with } r(e)\in (r(F)\backslash
s(F))\}],$$
$$\mbox{ and where } \ s((x,y))=x, \ \ r((x,y))=y \ \mbox{ for any }(x,y)\in E_{F}^{1}.$$
Note that, since $F$ is finite, the graph $E_F$ is finite (regardless of the size of $E$).}
\end{definition}

Remarks:  1.  It is conventional to define $s(v)=v$ for each vertex $v$ in $E$.  Because of that, the expression in rectangular brackets in Definition \ref{EsubFdefinition} for $E_F^1$ is redundant.  However, we choose to keep this expression in the definition, as it makes the correspondence
 between $E_F^1$ and the set $G^1$ in the proof of Proposition \ref{KeyLemma}  more transparent.

 2.  While the construction presented in Definition \ref{EsubFdefinition} is similar to that given in \cite[Definition 1.1]{RS}, there are indeed some significant differences.  For instance, the construction of \cite[Definition 1.1]{RS} requires that the graph $E$ has no sinks, while the construction presented here has no such stipulation.   Additionally, even in situations where $E$ is a graph with no sinks and $F$ is a finite subset of $E^1$, the two constructions can in fact yield different corresponding graphs $E_F$.  However, the underlying goal of each of the two constructions is the same, namely, to produce a subalgebra of a graph algebra which is isomorphic to the graph algebra of a finite graph.

\begin{example}\label{EsubFexample}
{\rm  For clarity, we provide an example of the graph $E_F$ constructed in the previous definition.  Let $E$ be the ``rose with $n$-petals" graph
$$E \ = \ 
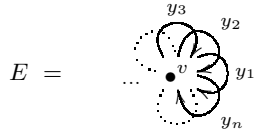
\xymatrix{ & {\bullet^v} \ar@(ur,dr) ^{y_1} \ar@(u,r) ^{y_2} \ar@(ul,ur) ^{y_3} \ar@{.} @(l,u) \ar@{.} @(dr,dl)
\ar@(r,d) ^{y_n} \ar@{}[l] ^{\ldots} }$$
  Let $F=\{y_1\}$.   Then $E_F^0 = \{y_1\} \cup \{v\}$, and $E_F^1 = \{(y_1,y_1),(y_1,v)\}$.  Pictorially, $E_F$ is given by
$$ E_F \ = \ \xymatrix{{\bullet}^{y_1} \ar@(ul,dl) _{(y_1,y_1)} \ar[r] ^ {(y_1,v)} & {\bullet}^v}$$

This  example indicates that various properties of the graph $E$ need not pass to the graph $E_F$.  For instance, $E$ is cofinal, while $E_F$ is not.  In particular, $L_K(E)$ is a simple algebra, while $L_K(E_F)$ is not.   (See \cite{AA3} for a more complete discussion.)}

\end{example}

Our interest in the construction given in Definition \ref{EsubFdefinition} can be generally described as follows.  We seek to
place each finite set of elements taken from the Leavitt path algebra $L_K(E)$ inside a subalgebra of $L_K(E)$ which  possesses
certain 'finiteness' properties.     In case $E$ is row-finite, by \cite[Lemma 3.2]{AMP} we can realize $L_K(E)$ as the direct union of subalgebras
of the form $L_K(E_i)$ where each $E_i$ is a finite, complete subgraph of $E$.  In the general case, however, we need not have such
a description of $L_K(E)$.   For instance, if $\aleph$ is an infinite cardinal, and ${\rm Clock}(\aleph)$ denotes the 'infinite clock' graph
$$ 
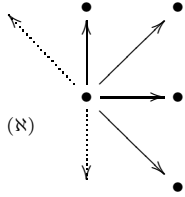
\xymatrix{ & {\bullet} & {\bullet} \\  & {\bullet} \ar@{.>}[ul] \ar[u] \ar[ur] \ar[r] \ar[dr] \ar@{.>}[d]
\ar@{}[dl] _{(\aleph)} & {\bullet} \\ &  & {\bullet}}$$
having $\aleph$ edges, then there are {\it no} nontrivial finite complete subgraphs of ${\rm Clock}(\aleph)$.

\begin{example}\label{EsubFexampleInfiniteClock}
{\rm  It will be instructive to consider the $E_F$ construction of Definition \ref{EsubFdefinition} within the infinite clock graph $E={\rm Clock}(\aleph)$.    So let $v$ denote the center vertex, let $f$ denote one of the edges, and let $w$ denote $r(f)$.  Let $F = \{f\}$.    Then $E_F^0 = \{f\} \cup \{w\}$, while $E_F^1 = \{(f,w)\}$.  Thus $E_F$ is the graph
$$ E_F \ = \ \xymatrix{{\bullet}^{f}
  \ar[r] ^ {(f,w)} & {\bullet}^w}$$
with two vertices, and one edge connecting them.   In particular, $L(E_F)\cong {\rm M}_2(K)$.}
\end{example}

Although in general $E_F$ need not be a subgraph of $E$ (indeed, as seen in Example~\ref{EsubFexample}, $E_F$ may contain more vertices than does $E$), there is an  important relationship between the Leavitt path algebras $L_K(E_F)$ and $L_K(E)$, as we now show.

\begin{proposition}\label{KeyLemma}    Let $F$ be a finite set of edges in a graph $E$. Then there is an algebra
homomorphism $\theta : L_{K}(E_{F})\rightarrow L_{K}(E)$ having the properties:

(1) $F \cup F^*  \subseteq {\rm Im}(\theta)$.

  (2) If $w\in r(F)$, then $w\in {\rm Im}(\theta)$.

  (3)  If $w\in E^0$ has $s_E^{-1}(w) \neq \emptyset$ and $s_E^{-1}(w) \subseteq F$, then $w \in {\rm Im}(\theta)$.

\end{proposition}

\begin{proof}
We define subsets $G^0$ and $G^1$ of $L_{K}(E)$ as follows.

\begin{eqnarray*}
&G^0&=\{ee^{\ast } \mid e\in F\} \\
&&  \ \ \ \ \cup \ \ \{v-\sum_{f\in F,s(f)=v}ff^{\ast } \mid v\in r(F)\cap s(F)\cap s(E^{1}\backslash F)\}\\
&&  \ \ \ \ \cup  \ \ \{v \mid v\in r(F)\backslash s(F)\} \\
\end{eqnarray*}

and

\begin{eqnarray*}
&G^1&=\{ eff^{\ast}  \mid e,f\in F,s(f)=r(e)\} \\
&& \ \ \ \ \cup \ \  \{e-
\sum_{f\in F,s(f)=r(e)}eff^{\ast } \mid r(e)\in r(F)\cap s(F)\cap s(E^{1}\backslash F)\}\\
&& \ \ \ \ \cup \ \
 \{e\in F \mid r(e)\in r(F)\backslash s(F)\}. \\
\end{eqnarray*}

 \medskip

 We define $\theta: L_K(E_F)\rightarrow L_K(E)$ as follows.

 \smallskip

 There are three different types of vertices in $E_F$.   If $w\in E_F^0$ has form $w = e\in F$, then define
 $$\theta(w) = ee^*.$$
  If $w\in E_F^0$ has form $w=v$ with $v\in r(F)\cap s(F)\cap s(E^{1}\backslash F)$, then define $$\theta(w) = v-\sum_{f\in F,s(f)=v}ff^{\ast }.$$
    If $w\in E_F^0$ has form $w=v$ with $v\in r(F)\backslash s(F)$, then define
    $$\theta(w)= w.$$

       Note that in each case we have $\theta(w)\in G^0$.

 \smallskip
 There are three different types of edges in $E_F$.
   If $h \in E_F^1$ has form $h = (e,f)$ with $f\in F$, then define
   $$\theta(h) = eff^*.$$
    If $h\in E_F^1$ has form $h=(e,r(e))$ with $r(e)\in r(F)\cap s(F)\cap s(E^{1}\backslash F)$, then define
    $$\theta(h)= e-
\sum_{f\in F,s(f)=r(e)}eff^{\ast }.$$
  If $h\in E_F^1$ has form $h=(e,r(e))$  with $r(e)\in r(F)\backslash s(F)$, then define
  $$\theta(h)=e.$$

      Note that in each case we have $\theta(h)\in G^1$.

\smallskip

For each $h\in E_F^1$ we define $\theta(h^*)=(\theta(h))^*$ in $L_K(E)$.

\medskip

It is now a long, straightforward check to verify that $\theta$ is compatible with the four types of relations which define  $L_K(E_F)$ (refer to Definition \ref{definition}).    As a representative example of the computations required here, we offer the following.   Let $w\in E^0_F$ have the form $w = e \in F$.  Then the (CK2) relation at $e$ in $L_K(E_F)$ is the equation
$$\sum_{g\in E_F^1, s(g)=e}gg^* = e.$$
   But $s(g)=e$ in $E_F^1$ means $g=(e,f)$ where either $f\in F$ has $s(f)=r(e)$, or $g = (e,r(e))$ with $r(e)\in r(F)\cap s(F)\cap s(E^{1}\backslash F)$, or $g = (e,r(e))$ with $r(e)\in r(F)\backslash s(F)$.  So the (CK2) relation at $e$ in $L_K(E_F)$ takes the form

\begin{eqnarray*}
&e&= \sum_{f\in F, s(f)=r(e)}(e,f)(e,f)^* \ + \sum_{w\in r(F)\cap s(F)\cap s(E^{1}\backslash F), w=r(e)}(e,w)(e,w)^* \\
&&  \ \ + \sum_{w\in r(F)\backslash s(F), w=r(e)}(e,w)(e,w)^*. \\
\end{eqnarray*}


Note that empty sums are interpreted as $0$.  Also, the final two summation expressions are in fact either singletons or empty, depending on whether $r(e)\in r(F)\cap s(F)\cap s(E^{1}\backslash F)$ or $r(e)\in r(F)\backslash s(F)$.

   We must show that the corresponding equation under $\theta$ holds in $L_K(E)$.  In other words, we must show

\begin{eqnarray*}
&ee^*&= \sum_{f\in F, s(f)=r(e)}(eff^*)(eff^*)^* \ \\
&&  \ \ +  \sum_{w\in r(F)\cap s(F)\cap s(E^{1}\backslash F), w=r(e)}[e - \sum_{f\in F,s(f)=r(e)} eff^*][e - \sum_{f\in F,s(f)=r(e)} eff^*]^* \\
&& \ \ + \sum_{w\in r(F)\backslash s(F), w=r(e)}ee^*. \\
\end{eqnarray*}


   There are two cases.  If $r(e)\in r(F)\backslash s(F)$, then this equation simply becomes $ee^* = ee^*$ and we are done.    On the other hand, if $r(e)\in s(F)$, then note the second 'sum' $\sum_{w\in r(F)\cap s(F)\cap s(E^{1}\backslash F), w=r(e)}[e - \sum_{f\in F,s(f)=r(e)} eff^*][e - \sum_{f\in F,s(f)=r(e)} eff^*]^*$ is in fact simply the single expression $[e - \sum_{f\in F,s(f)=r(e)} eff^*][e - \sum_{f\in F,s(f)=r(e)} eff^*]^*$.  So the right hand side is

\begin{eqnarray*}
&&\sum_{f\in F, s(f)=r(e)}(eff^*)ff^*e^* \ + \ [e - \sum_{f\in F,s(f)=r(e)} eff^*][e - \sum_{f\in F,s(f)=r(e)} eff^*]^*  \\
&=& \sum_{f\in F, s(f)=r(e)}eff^*e^* \ + \ [ee^* - \sum_{f\in F,s(f)=r(e)} eff^*e^*] \ \ \ \ \ \  \  \mbox{ (by (CK1) and (CK2)) }\\
&=& ee^* . \\
\end{eqnarray*}




In a similar manner one can verify the compatibility of $\theta$ with all the remaining relations which define $L_K(E_F)$.  Thus we conclude that  $\theta$ extends to an algebra homomorphism
$$\theta: L_K(E_F)\rightarrow L_K(E).$$

\medskip

By definition we have ${\rm Im}(\theta)$ is the subalgebra of $L_K(E)$ generated by $G^0,G^1,(G^1)^*$.  It will be helpful later to note that for each $x\in G^1 \cup (G^1)^*$ there exist $y,y' \in G^0$ for which $yxy' = x$.  In particular, if an element $z\in L_K(E)$ is orthogonal to every element of $G^0$, then necessarily $z$ is orthogonal to every element in ${\rm Im}(\theta)$.

\medskip

We are now in position to verify the three claimed properties of ${\rm Im}(\theta)$. For (1), we show that every $f\in F$ is contained in ${\rm Im}(\theta)$.
Suppose first that $f\in F$ has $s_E^{-1}(r(f))\subseteq F$.   If $s_E^{-1}(r(f)) = \emptyset$, then $r(f)\in r(F)\backslash s(F)$ vacuously, so by definition $f\in G^1 \subseteq {\rm Im}(\theta)$.   On the other hand, if $s_E^{-1}(r(f)) \neq \emptyset$, then we have $fgg^* \in G^1$ for all $g\in F$ having $r(f)=s(g)$.  But then by hypothesis this is the same as the collection of $g\in E^1$ having $r(f)=s(g)$.   Thus we have $\{fgg^* \mid g\in E^1, s(g)=r(f)\}\subseteq {\rm Im}(\theta)$, so that in particular ${\rm Im}(\theta)$ contains  $\sum _{g\in E^1, s(g)=r(f)}fgg^* = f \cdot \sum _{g\in E^1, s(g)=r(f)}gg^* = f\cdot r(f)=f$.

\smallskip

On the other hand, suppose that $f\in F$ has the property that $s_E^{-1}(r(f)) \nsubseteq F$.
Then there are two possibilities.  In the first case,  $s_E^{-1}(r(f)) \neq \emptyset$ and  $s_E^{-1}(r(f)) \cap F = \emptyset$.  (In other words, there are edges in $E$ which are emitted from $r(f)$, but none of these edges are in $F$.)   But then $r(f) \notin s(F)$, so that $f\in G^1$ by definition, so that $f\in {\rm Im}(\theta)$.    In the second case, suppose  $s_E^{-1}(r(f)) \cap F \neq \emptyset$.  Then either we have $s_E^{-1}(r(f)) \subseteq F$ (in which case we are done by the previous paragraph), or we have $r(f)\in r(F)\cap s(F)\cap s(E^{1}\backslash F)$.  In this situation we have $fgg^*\in G^1\subseteq {\rm Im}(\theta)$ for all $g\in F$ having $s(g)=r(f)$, so we in particular have $\sum _{g\in F, s(g)=r(f)}fgg^*$ in ${\rm Im}(\theta)$ as well.  But by definition we also have the element $f-\sum _{g\in F, s(g)=r(f)}fgg^*$ in $G^1\subseteq {\rm Im}(\theta)$.  Then
$$f = (\sum _{g\in F, s(g)=r(f)}fgg^*) \ + \ (f-\sum _{g\in F, s(g)=r(f)}fgg^*) \in {\rm Im}(\theta).$$
Thus we conclude that $F\subseteq {\rm Im}(\theta)$.  But for each $x\in {\rm Im}(\theta)$ we have $x^* \in {\rm Im}(\theta)$ by definition.  Thus $F\cup F^* \subseteq {\rm Im}(\theta)$, thereby establishing (1).

\medskip

   In particular, if $w = r(f)$ for $f\in F$, then $w = f^*f \in {\rm Im}(\theta)$, yielding (2).  For (3), suppose $s^{-1}(w) \neq \emptyset$, and
$s^{-1}(w)\subseteq F$.  Then each $ff^*$ for $f\in E^1$ having $s(f)=w$  is in $G^0$, so that $ \sum_{f\in E^1, s(f)=w}ff^*$ is in ${\rm Im}(\theta)$.  But this last sum is precisely $w$ by (CK2).


\end{proof}

We remark here that for $\theta$ as given in Proposition \ref{KeyLemma}, $\theta(w)\neq 0$ for all three possible types of $w\in E_F^0$.  (That $\theta(w)\neq 0$ in case $w\in r(F)\cap s(F)\cap s(E^{1}\backslash F)$ hinges on the fact that there exists $g\in E^{1}\backslash F$ having $s(g)=w$.)   This in turn will allow us to conclude, in certain situations (including the situation where $E$ is acyclic), that $\theta$ is in fact a monomorphism.  (See e.g. \cite{AALP}.)   However, we will not utilize this additional property of  $\theta$ in the sequel.

\medskip

With Proposition \ref{KeyLemma} in hand, we now construct the subalgebras of $L_K(E)$ which play the central role in our main result, Theorem \ref{VNRegular}.

\medskip

\subsection*{{\bf The Subalgebra Construction}}

Let $E$ be any graph, $K$ any field, and $\{a_1,a_2,...,a_{\ell}\}$ any finite subset of nonzero elements of $L_K(E)$.   For each $1 \leq r \leq \ell$ write
$$a_r = k_{c_1}v_{c_1} + k_{c_2}v_{c_2} + ... + k_{c_{j(r)}}v_{c_{j(r)}} + \sum_{i=1}^{t(r)} k_{r_i}p_{r_i}q^*_{r_i}$$
where each $k_j$ is a nonzero element of $K$, and, for each $1\leq i \leq t(r)$, at least one of $p_{r_i}$ or $q_{r_i}$ has length at least 1.   (That such a representation for each $a_r$ exists follows from properties of $L_K(E)$ mentioned previously.)
Let  $F$ denote the (necessarily finite) set of those edges in $E$ which appear in the representation
of some $p_{r_i}$ or $q_{r_i}$, $1\leq r_i \leq t(r)$, $1\leq r \leq \ell$.
Now consider the set
$$S = \{v_{c_1}, v_{c_2}, ... , v_{c_{j(r)}} \mid 1\leq r \leq \ell \}$$
 of vertices which appear in the displayed description of $a_r$ for some $1\leq r \leq \ell$.
We partition $S$ into subsets as follows:
    $$S_1 = S \cap r(F),$$
and, for  the remaining vertices $T = S\backslash S_1$, we define
       $$S_2 = \{v\in T \mid s_E^{-1}(v)\subseteq F \mbox{ and } s_E^{-1}(v) \neq \emptyset \}$$
$$S_3 = \{v\in T \mid s_E^{-1}(v) \cap F = \emptyset \}$$
$$S_4 = \{v\in T \mid s_E^{-1}(v) \cap F \neq \emptyset \mbox{ and } s_E^{-1}(v) \cap (E^1 \backslash F) \neq \emptyset \}.$$
Let $E_F$ be the  graph as constructed in Definition \ref{EsubFdefinition} corresponding to this set $F$, and let $\theta: L_K(E_F) \rightarrow L_K(E)$ be the homomorphism described in Proposition \ref{KeyLemma}.

\begin{definition}\label{B(S)definition}
{\rm Let $E$ be any graph, $K$ any field, and $\{a_1,a_2,...,a_{\ell}\}$ any finite subset of nonzero elements of $L_K(E)$.  Consider the notation presented in The Subalgebra Construction.   We define  $B(a_1,a_2,...,a_{\ell})$ to be the $K$-subalgebra of $L_K(E)$ generated  by the set ${\rm Im}(\theta) \cup S_3 \cup S_4$.  That is,
$$B(a_1,a_2,...,a_{\ell}) = <{\rm Im}(\theta), S_3, S_4>.$$
}
\end{definition}

\begin{proposition}\label{B(S)}
Let $E$ be any graph, $K$ any field, and $\{a_1,a_2,...,a_{\ell}\}$ any finite subset of nonzero elements of $L_K(E)$.   Let $F$ denote the subset  of $E^1$ presented in The Subalgebra Construction.  For $w\in S_4$ let $u_w$ denote the element $w - \sum_{f\in F, s(f)=w} ff^*$ of $L_K(E)$.  Then

(1)  $\{a_1,a_2,...,a_{\ell}\} \subseteq B(a_1,a_2,...,a_{\ell})$.

(2)  $B(a_1,a_2,...,a_{\ell})={\rm Im}(\theta) \oplus (\oplus_{v_i\in S_3} Kv_i ) \oplus (\oplus_{w_j \in S_4} Ku_{w_j} ).$

(3)  The collection $\{B(S)\mid S\subseteq L_K(E), S \mbox{ finite}\}$ is an upward directed set of subalgebras of $L_K(E)$.

(4)  $L_K(E) = \mathop{\oalign{\hfill $\limind $\hfill\cr}}_{{\{S\subseteq L_K(E),S \ finite \}}} B(S)$.


\end{proposition}

\medskip

\begin{proof}

(1)  By Proposition \ref{KeyLemma} we have that $F \cup F^* \cup S_1 \cup S_2 \subseteq  B(a_1,a_2,...,a_{\ell})$.  Since  $S_3\cup S_4 \subseteq  B(a_1,a_2,...,a_{\ell})$ by construction, (1) follows.

\medskip

(2)  Since the $\{v_i\}$ and $\{u_{w_j}\}$ are pairwise orthogonal idempotents we immediately get that $\sum_{v_i\in S_3} Kv_i  = \oplus_{v_i\in S_3} Kv_i$, and that $ \sum_{w_j \in S_4} Ku_{w_j} = \oplus_{w_j \in S_4} Ku_{w_j}.$   We now establish that the three indicated summands are mutually orthogonal, which will establish that the sum $
{\rm Im}(\theta) + (\oplus_{v_i\in S_3} Kv_i ) + (\oplus_{w_j \in S_4} Ku_{w_j})$  is direct.
  Let $v\in S_3$.  Then by definition $v$ is neither the source vertex nor range vertex for any of the elements in $F$.  In particular, the one dimensional subalgebra $Kv$ of $L_K(E)$ clearly annihilates all the elements of $G^0$; as noted previously, this suffices to yield that $Kv$ indeed annihilates ${\rm Im}(\theta)$.
 But for any $w\in S_4$ we have that $u_w = w - \sum_{f\in F, s(f)=w} ff^*$  is orthogonal in $L_K(E)$ to $v$, since $S_3\cap S_4 = \emptyset$.

So we have shown that $\oplus_{v_i\in S_3} Kv_i \ \cap \ [{\rm Im}(\theta) + \oplus_{w_j \in S_4} Ku_{w_j}] = \{0\}$.  Thus we need only show that ${\rm Im}(\theta) \cap Ku_w = \{0\}$ for all $w\in S_4$, which we establish by showing $Ku_w\cdot {\rm Im}(\theta) = {\rm Im}(\theta)\cdot Ku_w = \{0\}$ and using the fact that each $u_w$ is idempotent.   Choose any such $w$.    Since by definition $w \notin r(F)$ we have that $u_w = w - \sum_{f\in F, s(f)=w} ff^*$ is orthogonal in $L_K(E)$ to elements of $G^0$ of the form $w' - \sum_{g\in F, s(g)=w'} gg^*$ for $w'\in r(F)\cap s(F)\cap s(E^1\backslash F)$.  Similarly, $u_w$ is orthogonal in $L_K(E)$ to $r(e)$ for each $e\in F$.   Now suppose $ee^* \in G^0$ with $e\in F$.   If $s(e)\neq w$ then $w - \sum_{f\in F, s(f)=w} ff^*$ is clearly orthogonal to $ee^*$.  On the other hand, if $s(e) = w$ then $u_w\cdot ee^* = (w - \sum_{f\in F, s(f)=w} ff^*)\cdot ee^* = wee^* - \sum_{f\in F, s(f)=w} ff^*ee^* =  ee^* - ee^* = 0$, with the final simplification occurring because $ff^*ee^* = ee^*$ for $e=f$, and $ff^*ee^*=0$ otherwise by (CK1).  Similarly, we have $ee^* \cdot u_w = 0$.  We conclude that the indicated sums are direct.

We now show that the direct sum in fact equals $B(a_1,a_2,...,a_{\ell})$. By construction,  it suffices to show that $u_{w}\in B(a_1,a_2,...,a_{\ell})$ for all $w\in S_4$, and that $w \in {\rm Im}(\theta) \oplus (\oplus_{v_i\in S_3} Kv_i ) \oplus (\oplus_{w_j \in S_4} Ku_{w_j} )$ for all $w\in S_4$.  But each of these inclusions follow directly by noting that,
for each $w \in S_4$ and each $f\in F$ having $s(f)=w$, we have $ff^* \in {\rm Im}(\theta)$ by definition.

\medskip

(3)  As noted previously, $E_F$ is a finite graph for each finite subset $F$ of $E^1$. In particular, $L_K(E_F)$ is a finitely generated $K$-algebra (with generating set $E_F^0 \cup E_F^1 \cup (E_F^1)^*$).  This in turn implies that ${\rm Im}(\theta)$, and hence $B(S)$, is a finitely generated $K$-algebra for each finite set $S$ of $L_K(E)$.   In particular, if $S_1$ and $S_2$ are finite subsets of $L_K(E)$, we let $T_1$ (resp. $T_2$) denote a finite set of generators of $B(S_1)$ (resp. $B(S_2)$).  If $T=T_1\cup T_2$, then it is clear by construction that $B(S_1) \cup B(S_2) \subseteq B(T)$.

\smallskip

(4) now follows immediately from (1) and (3).

\end{proof}

As noted previously, various properties of the graph $E$ need not pass to the graph $E_F$.    However,

\begin{lemma}\label{acyclic}
Let $E$ be any acyclic graph, and $F$ any finite subset of $E^1$.  Then $E_F$ is acyclic.
\end{lemma}

  \begin{proof}
  By contradiction, we show that the existence of a closed path in $E_F$ necessarily yields a closed path in $E$. By
definition, a closed path in $E_F$ is of the form $(e_1,e_2),(e_2,e_3), ... ,(e_n,e_1)$ where $(e_i,e_{i+1})\in E_F^1$. Now it is straightforward to
show that the indicated sequence of edges in $E_F$ yields a sequence $e_1, e_2, ..., e_n$ in $E^1$ having the desired property.
\end{proof}

We recall now some ideas which play central roles in our main result.  For additional information about these concepts, see for example \cite{G2}, \cite{K}, and  \cite{NZ}.

\begin{definition}  {\rm  Let $R$ be a (not necessarily unital) ring.

  (1)  $R$ is called  {\it von Neumann regular} in case for every $x\in R$ there exists $y\in R$ such that
 $x=xyx$.

 (2)  $R$ is called $\pi$-{\it regular} in case for every $x\in R$ there exists $y\in R$ and $n \in \N$ for which $x^n = x^nyx^n$.

  (3) $R$ is called \textit{left} (resp. \textit{right}) $\pi $-
\textit{regular} if for each $a\in R$ there exists $n\in {\mathbb{N}}$ and $b\in R$ such that $a^{n}=ba^{n+1}$ (resp. $a^{n}=a^{n+1}b$).   (For rings with
local units, this is equivalent to saying that the descending chain of left
ideals $Ra\supseteq Ra^{2}\supseteq ...\supseteq Ra^{k}\supseteq ...$ \ (resp. right ideals $aR\supseteq a^{2}R\supseteq ...\supseteq a^{k}R\supseteq ...$)
becomes stationary after finitely many terms.)

(4)  $R$ is called \textit{strongly }$\pi -$\textit{regular} if its both
left and right $\pi $-regular.
}

\end{definition}

Clearly any von Neumann regular ring is $\pi$-regular.  Conversely,  the ring $\Z / 4\Z$ provides an easy example of a ring which is $\pi$-regular but not von Neumann regular (since $\bar{2}$ has no von Neumann regular inverse).

 By \cite[Lemma 6]{CY}, if $R$ is a unital strongly $\pi$-regular ring then
for every
element $a\in R$ there is a positive integer $n$ and an element $x\in R$
such that $ax=xa$ and $a^{n+1}x=a^{n}=xa^{n+1}$.  (We will show below that this result holds for rings with local units as well.)  From this property it is then easy to show that if $R$ is strongly $\pi$-regular, then $R$ is $\pi$-regular.
Conversely, the ring $R={\rm End}_K(V)$ of all linear transformations of an infinite
dimensional vector space $V$ over a field $K$ provides an example of a ring which is $\pi$-regular (in fact, von Neumann regular), but  not strongly $\pi$-regular.   (Indeed, if $\alpha: V\rightarrow V$ is the shift transformation given by $\alpha (K_{1})=0$ and
$\alpha (K_{i+1})=K_{i}$ for all $i>1$, then for any $n$, $\ker \alpha
^{n}=\oplus _{i=1}^{n}K_{i}$ and so $\alpha ^{n}\neq \beta \alpha ^{n+1}$
for any $n.$)

\begin{lemma}\label{localstronglypiregular}
 Let R be
a ring with local units.  Then
R is  strongly $\pi$-regular ring if and only if for every nonzero idempotent $v$ of $R$, the subring $vRv$ is strongly $\pi$-regular.
\end{lemma}

\begin{proof}

Assume R is strongly $\pi$-regular.  Pick  $a \in vRv$.    By hypothesis there exists $b \in R$
with $a^n=a^{n+1}b$,   and there exists  $c \in R$ with $a^m = ca^{m+1}$.        But
$vav = v$,  so  $va^nv = a^n$ and $va^{n+1}v = a^{n+1}$.  Thus, multiplying both sides of the equation  $a^n=a^{n+1}b$ by $v$, we get $a^n = va^{n+1}vbv = a^{n+1}vbv$.   Since $vbv \in vRv$ we
have shown that  $a^n=a^{n+1}b'$  for some $b' \in vRv$.    A similar computation
yields that $a^m = c' a^{m+1}$  inside $vRv$.

Conversely, pick $a \in R$.   Then $a \in vRv$ for some idempotent $v$ by definition of set
of local units.   So there exist $b$  and $c$ in $vRv$, and hence in $R$, with the appropriate properties.
\end{proof}

Although the properties von Neumann regular, strongly $\pi$-regular, and $\pi$-regular are in general not equivalent, as one consequence of Theorem \ref{VNRegular} we conclude that these properties are indeed equivalent in the context of Leavitt path algebras.

\smallskip

We are now in position to establish our main result.

\begin{theorem}\label{VNRegular}   Let $E$ be an arbitrary graph, and let $K$ be any field.  The following are equivalent.

 (1)  $L_K(E)$ is von Neumann regular.

 (2)  $L_K(E)$ is $\pi$-regular.

 (3)  $E$ is acyclic.

 (4)  $L_K(E)$ is locally $K$-matricial; that is, $L_K(E)$ is the direct union of subrings, each of which is isomorphic to a finite direct sum of finite matrix rings over $K$.

 (5)  $L_K(E)$ is strongly $\pi$-regular.

\end{theorem}

\begin{proof}

 (1) $\Rightarrow$ (2) is immediate.

\medskip

(2) $\Rightarrow$ (3).    By contradiction, suppose $c $ is a cycle in $E$, and let $v=s(c)=r(c)$.   We show that  $v+c$ has no $\pi$-regular inverse in $L_K(E)$.

 Let $\gamma$ denote $v+c$.  Suppose there exist  $\beta \in L_{K}(E)$, $n\in \N$
such that $\gamma ^{n}\beta \gamma ^{n}=\gamma ^{n}$. Note that, since
$\gamma v=\gamma =v\gamma $, we have $\gamma ^{n}v\beta v\gamma ^{n}=\gamma
^{n}$. Then $\alpha =v\beta v$ satisfies $\gamma ^{n}\alpha \gamma
^{n}=\gamma ^{n}$ and $v\alpha v=\alpha $. Moreover, $\alpha v=v\alpha
=\alpha $.

Write $\alpha $ as a graded sum $\alpha =\sum\limits_{i=M}^{N}a_{i}$ where $a_{M}\neq 0$, $a_N \neq 0$, $\deg (a_{i})=i$ for all nonzero $a_i$ having
 $M\leq i\leq N$, and $a_{i}=0$ if $i>N$ or $i<M$. Since $\deg (v)=0,$
the equation $\alpha v=v\alpha $ implies that $a_{i}v=va_{i}=a_{i}$ for all $%
i$. \ Now expanding the equation $\gamma ^{n}\alpha \gamma ^{n}=\gamma ^{n}$%
, we obtain
\[
(v+\sum\limits_{k=1}^{n}
\binom{n}{k}
c^{k})(\sum\limits_{i=M}^{N}a_{i})(v+%
\sum\limits_{k=1}^{n}\binom{n}{k}c^{k})=v+\sum\limits_{k=1}^{n}\binom{n}{k}%
c^{k}.
\]
Equating the lowest degree terms on both sides, we get $va_{M}v=v$, so that $%
a_{M}=v$. Since $\deg (v)=0,$ we conclude that $M=0$, and that $a_{0}=v$. \ Thus $%
\alpha =\sum\limits_{i=0}^{N}a_{i}$. Let $\deg (c)=s>0$. Now every term
other than the first on the right hand side has degree $ks$ for some
positive integer $k\leq n$, and so equating the corresponding graded
components on both sides, we conclude that $a_{i}=0$ if $i$ is not a
multiple of $s.$

We establish by induction  that $a_{ks}=f_{k}(c)$ for each $k \in \N$, where $f_k(c)$ is a polynomial in $c$ with
integer coefficients. For $k=1$, by equating the degree $s$ components on both
sides we obtain
$$
va_{s}v+\binom{n}{1}ca_{0}+a_{0}\binom{n}{1}c=\binom{n}{1}c.
$$
This implies that $a_{s}=-\binom{n}{1}c$, an integral polynomial in $c$.
Now suppose $t>1$, and suppose $a_{ks}=f_{k}(c)$, an integral
polynomial in $c$ for all $1\leq k<t.$  We expand the previously displayed equation, and equate the
degree $ts$ terms of both sides.  This yields

\begin{eqnarray*}
&a_{ts}&+\binom{n}{1}c[a_{(t-1)s}+a_{(t-2)s}\binom{n}{1}c+a_{(t-3)s}\binom{n}{%
2}c^{2}+ ... +a_{0}\binom{n}{t-1}c^{t-1}] \\
&&+\binom{n}{2}c^{2}[a_{(t-2)s}+a_{(t-3)s}\binom{n}{1}c+ ... +a_{0}\binom{n}{%
t-2}c^{t-2}] \\
&& +\binom{n}{3}c^{3}[a_{(t-3)s}+a_{(t-4)s}\binom{n}{1}c+ ... +a_{0}\binom{n}{%
t-3}c^{t-3}] \\
&& + \ ... \ + \ \binom{n}{t}c^{t} a_{0} \\
&=&\binom{n}{t}c^{t}
\end{eqnarray*}%

Substituting for $a_{s},...,a_{(t-1)s}$ as allowed by the induction hypothesis and solving for $a_{ts}$, we obtain $%
a_{ts}=f_{t}(c)$, a polynomial in $c$ with integer coefficients.

In particular, we conclude that every homogeneous component $a_i$ of $\alpha$ commutes with $c$ in $L_K(E)$.  This yields
that $c\alpha = \alpha c$.   But  then the equation $(v+c)^{n}\alpha
(v+c)^{n}=(v+c)^{n}$ becomes
$$
\alpha (v+c)^{2n}=(v+c)^{n}.$$

But this is impossible, as follows.   Since each $a_i$ is a polynomial in $c$ with integer coefficients, we have $a_i c^r \neq 0$
for all $r\in \N$.    Let $i$ be maximal with the property that $a_i(v+c)^{2n} \neq 0$.  (Such $i$ exists, since
$a_0 = v$ has this property.)   Then the left hand side contains terms of degree $2sn + i$ (namely, $a_ic^{2n}$), while the maximum
degree of terms on the right hand side is $ns$.

\medskip

(3) $\Rightarrow$ (4).    We assume $E$ is acyclic.    Let $\{B(S)\mid S\subseteq L_K(E), S \mbox{ finite}\}$ be the collection of subalgebras of $L_K(E)$ indicated in Proposition \ref{B(S)}(3).   By Proposition \ref{B(S)}(4), it suffices to show that each such  $B(S)$ is of the indicated form.   But by Proposition \ref{B(S)}(2), $B(S)={\rm Im}(\theta) \oplus (\oplus_{v_i\in S_3} Kv_i ) \oplus (\oplus_{w_j \in S_4} Ku_{w_j} ).$  Since terms appearing in the second and third summands are clearly isomorphic as algebras to $K\cong {\rm M}_1(K)$, it suffices to show that ${\rm Im}(\theta)$ is isomorphic to a finite direct sum of finite matrix rings over $K$.   Since $E$ is acyclic, by Lemma \ref{acyclic} we have that $E_F$ is acyclic.  But $E_F$ is always finite by definition, so we have by \cite[Proposition 3.5]{AAS} that $L(E_F)\cong  \oplus_{i=1}^{\ell} {\rm M}_{m_i}(K)$ for some $m_1,...,m_{\ell}$ in $\N$.  Since each ${\rm M}_{m_i}(K)$ is a simple ring, we have that any homomorphic image of $L_K(E_F)$ must have this same form.  So
we get that  ${\rm Im}(\theta) \cong \oplus_{i=1}^{L} {\rm M}_{m_i}(K)$ for some $m_1,...,m_{L}$ in $\N$, and we are done.   (As remarked previously, since $\theta$ is in fact an isomorphism we have $L = \ell$.)

   \medskip

(4) $\Rightarrow$ (1).  It is well known that any algebra of the form $\oplus_{i=1}^{\ell} {\rm M}_{m_i}(K)$ is von Neumann regular.   But every element of $L_K(E)$ is contained in a subalgebra of $L_K(E)$ of this form, so that every element of $L_K(E)$ thereby has a von Neumann regular inverse.

\medskip

(4) $\Rightarrow$ (5).  Suppose $L_{K}(E)$ is locally $K$-matricial. So every element $a\in L_{K}(E)$
is contained in a subring $S\cong \oplus_{i=1}^{\ell} {\rm M}_{m_i}(K)$. As any such $S$ is a unital left (resp. right) artinian ring, there
is a $b\in S$ and a positive integer $n$ such $a^{n}=ba^{n+1}$ (resp. $
a^{n}=a^{n+1}b$).

\medskip

(5) $\Rightarrow$ (2)  By Lemma \ref{localstronglypiregular} we have that each $a\in L_K(E)$ is contained in a strongly $\pi$-regular unital subring of the form $vL_K(E)v$ for some $v=v^2\in L_K(E)$.   Then by \cite[Lemma 6]{CY}  there is a positive integer $n$ and an element $x\in vL_K(E)v$
such that $ax=xa$ and $a^{n+1}x=a^{n}=xa^{n+1}$.  Now iterating the substitution $a^n = a^{n+1}x = aa^nx=a(a^{n+1}x)x=a^{n+2}x^2$ we get $a^{n}=a^{2n}x^n$, which using $ax=xa$ gives $a^n = a^nx^na^n$, which yields (2).

\end{proof}

We record the following consequence of Theorem \ref{VNRegular}, in part because it demonstrates the independence of our results from any cardinality restrictions or graph-theoretic restrictions (e.g. row-finiteness) on the graphs.

\begin{example}
{\rm Let $\aleph$ be any cardinal, and let ${\rm Clock}(\aleph)$ be the infinite clock graph having $\aleph$ edges.   Then for any field $K$, the Leavitt path algebra $L_K({\rm Clock}(\aleph))$ is von Neumann regular.  In addition, $L_K({\rm Clock}(\aleph))$ is locally $K$-matricial.
}
\end{example}

It is worth noting that the locally $K$-matricial nature of $L_K({\rm Clock}(\aleph))$ does {\it not} stem from a consideration of the finite complete subgraphs of ${\rm Clock}(\aleph)$, since as noted previously ${\rm Clock}(\aleph)$ contains {\it no} such nontrivial subgraphs.

\medskip

As a second consequence of  Theorem \ref{VNRegular},  we see that
 the ring $R = {\rm End}_K(V)$
of all linear transformations of an infinite dimensional vector space $V$
over a field $K$  cannot be represented as $L_K(E)$ for any
graph $E$, since $R$ is von Neumann regular but not strongly $\pi $-regular (as noted earlier).   Similarly,
let $V$ be a
vector space of uncountable dimension over a field $K$ and let $S$ be the (nonunital) subring of ${\rm End}_K(V)$ consisting of those linear transformations whose images are of at most countable
dimension.  Then $S$ is a von Neumann regular ring with local units.  However,
  $S$ is not strongly $\pi$-regular, so again invoking Theorem \ref{VNRegular} we have that $S$ cannot be
represented as the Leavitt path algebra of any graph $E$.



\medskip

We conclude this article by analyzing two additional ``regularity" properties of a ring.  We recall the definitions of some ring-theoretic terms.

 \begin{definition} {\rm  Let $R$ be a unital ring.

 (1) $R$ is called  {\it clean}  if each $a\in R$ is of the form $a=e+u$ where $e$ is an idempotent and $u$ is a (two-sided) unit. If in addition $aR\cap eR=0,$ we say $R$ is a {\it special
clean ring}. A clean ring $R$ is said to be {\it strongly clean} if in
the above definition we can choose $e$ and $u$ which commute.

(2)  $R$  is called {\it unit regular} in case for each $a\in R$ there exists a (two-sided) unit $u\in R$ such that $aua=a$.  In particular, every unit regular ring is von Neumann regular.
}
\end{definition}

Additional information about clean rings can be found in \cite{NZ}, while additional information about unit regular rings can be found in \cite{G2}.
The properties ``clean" and ``unit regular" are exemplified by matrix rings. Indeed if $R$ is the ring
 of $n\times n$ matrices over a field, then $R$ is both unit regular \cite[page 38]{G2} and strongly clean \cite[Theorem 4.1]{NZ}.  By \cite[Theorem 1]{CK}, a unital ring $R$ is unit regular if and only if $R$ is a special clean ring; in particular, any ring of the form ${\rm M}_n(K)$ for $K$ a field and $n\in \N$ is a special clean ring.

While the definitions of von Neumann regularity and $\pi$-regularity extend verbatim from unital rings to the nonunital case, the notions of clean and unit regularity require additional attention in the nonunital situation (since each definition refers to a unit in the given ring).   We now show how to naturally extend these latter two notions to rings with local units.

\begin{definition}
{\rm Let $R$ be a ring with local units.

 (1) $R$ is called {\it locally unit
regular} if for each $a\in R$ there is an idempotent $v\in R$ for which $a\in vRv$,
and elements
$u,u'\in vRv$ such that $uu^{\prime }=v=u^{\prime }u$, and $aua=a$.

 (2) $R$ is called {\it locally clean} if for each $a\in R$ there is an idempotent $v\in R$ for which $a\in vRv$, and elements $e,u,u'\in vRv$ such
 that $e$ is an idempotent, $uu' = v = u'u$, and $a=e+u$.
 }
\end{definition}

That the two notions given in the previous definition are natural generalizations of the corresponding notions for unital rings is established in the following.

\begin{lemma}  Let $R$ be a unital ring.

 (1)  $R$ is locally unit regular if and only if $R$ is unit regular.

  (2)  $R$ is locally clean if and only if $R$ is clean.

\end{lemma}

\begin{proof}

For (1), suppose $R$ is a ring with $1$ and is locally unit regular.
Let $a\in R$, and let $v,u,u'$ as given in the definition.  Then $w=u+(1-v)$ and $w^{\prime }=u^{\prime }+(1-v)$ satisfy $ww^{\prime
}=1=w^{\prime }w$ and $a=awa$. Hence $R$ is unit regular. The converse is clear with $v=1$.

Likewise, for (2), suppose $R$ is a ring with $1$ and is locally clean. Let $a\in R$, and write  $a=u+e$ as given in the definition.  Then $a=w+e^{\prime }$, where $e^{\prime }=e+(1-v)$ is an idempotent and $w=u-(1-v)$ is a two-sided unit in $R$ (since
with $w^{\prime }=u^{\prime }-(1-v)$, we have $ww^{\prime }=1=w^{\prime }w$).
Thus $R$ is clean.  As with (1), the converse follows with $v=1$.

\end{proof}

Our final result shows that for acyclic graphs $E$,  $L_{K}(E)$ possesses
the locally unit regular property, as well as  a property involving clean unital subrings.

\begin{theorem}
Let $E$ be an arbitrary graph, and let $K$ be any field. Then the following
conditions are equivalent:

(1) $E$ is is acyclic.

(2) $L_{K}(E)$ is locally unit regular.

(3) $L_{K}(E)$ is a direct limit of unital strongly clean rings, each of which is special.

\end{theorem}
\vspace{1pt}

\begin{proof}
(1) $\Leftrightarrow$ (2)  Suppose $E$ is acyclic. Then, by Theorem \ref{VNRegular}, $L_{K}(E)$ is a direct union of direct sums of matrix rings each of which, by \cite[page 38]{G2}, is unit regular. It is then clear that $L_{K}(E)$
is locally unit regular, where for each $a\in L_K(E)$ we use  for $v$ the identity element of the corresponding subring $B(S)$. Conversely, if $L_{K}(E)$ is locally unit regular, then it
is, in particular, von Neumann regular. So, by Theorem \ref{VNRegular}, $E$ is acyclic.

\smallskip

(2) $\Leftrightarrow$ (3) \ Suppose $L_{K}(E)$ is locally unit regular.
Since it is von Neumann regular, Theorem \ref{VNRegular} implies that $L_{K}(E)$ is a
directed union of direct sums of matrix rings $L_{i}$\ each of which, as
noted above, is a special clean ring which is, in addition, strongly
clean. On the other hand, if $L_{K}(E)$ is a directed union of special clean
rings $L_{i}$, then each $L_{i}$ is unit regular by \cite[Theorem 1]{CK},  and
so $L_{K}(E)$ is locally unit regular.  (Again for each $a\in L_K(E)$ we use  for $v$ the identity element of the corresponding subring $B(S)$.)
\end{proof}

A study of $L_{K}(E)$ for arbitrary graphs $E$ is presented by Goodearl in \cite{G}. Included in \cite{G} is a method to write $E$ as a
direct union of countable complete subgraphs.  We now show how this
approach together with the desingularization process yields an alternate proof
of the implication (3) $\Rightarrow $ (1) of Theorem \ref{VNRegular}. Our
aim in doing so is to contrast this approach with that of using Proposition \ref{KeyLemma}, which helps us to establish not only the von Neumann regularity of $L_{K}(E)$ for an acyclic graph $E$, but uncovers several internal properties
of such an $L_{K}(E)$ (e.g.,  locally matricial and locally unit regular). (Our approach also shows the coincidence of
von Neumann regularity with $\pi$-regularity and strong $\pi$-regularity
for Leavitt path algebras.) One may also note
that the desingularization approach as shown below does
not work for $\pi$-regular rings since $\pi$-regularity, unlike von
Neumann regularity, is not a Morita invariant (see e.g. \cite{R}). Proposition \ref{KeyLemma} poses no such restrictions, and provides additional structural insight into
these rings.

So suppose $E$ is acyclic.   By \cite[Proposition 2.7]{G}  $L_K(E) = \mathop{\oalign{\hfill $\limind $\hfill\cr}}_{\alpha \in A} L_K(E_{\alpha})$, with the limit taken over
the set $\{E_{\alpha} \mid \alpha \in A\}$ of countable complete subgraphs of $E$.  So in order to show that $L_K(E)$ is von Neumann regular, it suffices to show that each $L_K(E_{\alpha})$ is
von Neumann regular, since the direct limit of von Neumann regular rings is von Neumann regular. Since $E$ is
acyclic then necessarily so is each $E_{\alpha}$.

Since $E_{\alpha}$ is countable, we may form a desingularization $F_{\alpha}$ of $E_{\alpha}$.  (See e.g. \cite{AA3}.)   By construction, $F_{\alpha}$ is
row-finite.   Also, since desingularization preserves Morita equivalence, and von Neumann regularity is preserved by Morita
equivalence for rings with local units by \cite[Proposition 3.1]{AM}, it suffices to show that each $L_K(F_{\alpha})$ is von Neumann
regular. Since each $E_{\alpha}$ is acyclic, the desingularization construction shows that each $F_{\alpha}$ is acyclic as well.

But by \cite[Lemma 3.2]{AMP} $F_{\alpha}$ is the direct union of $G_{\beta}$ (the union taken over the set of finite complete subgraphs
of $F_{\alpha})$, and $L_K(F_{\alpha}) =
\mathop{\oalign{\hfill $\limind $\hfill\cr}}_{\beta \in B} L_K(G_{\beta})
 $.   Thus it suffices to show that each
$L_K(G_{\beta})$ is von Neumann regular.   Since $G_{\beta}$ is a subgraph of $F_{\alpha}$ we have that $G_{\beta}$ is acyclic.

So in the end, to establish that $L_K(E)$ is von Neumann regular,  it suffices to show that for any finite acyclic graph $G$ that $L_K(G)$ is von Neumann regular.   But by
\cite[Proposition 3.5]{AAS}  the Leavitt path algebra of a finite acyclic graph is isomorphic to a finite
direct sum of finite dimensional matrix rings over the ground field $K$, and such rings are well known to be von
Neumann regular (see e.g. \cite[Section 1]{G2}).

\medskip

We conclude this article by noting one more consequence of Theorem \ref{VNRegular} (we thank the referee for this suggestion).  The proof follows directly from the fact that von Neumann regularity is a Morita invariant for rings with local units.  We contrast this result with the aforementioned remark that, in general, the $\pi$-regularity property is not a Morita invariant.

\begin{corollary}
The property of $\pi$-regularity is a Morita invariant for Leavitt path algebras; that is, if $E$ and $F$ are graphs with $L_K(E)$ Morita equivalent to $L_K(F)$, then $L_K(E)$ is $\pi$-regular if and only if $L_K(F)$ is $\pi$-regular, and in this case $E$ and $F$ are both acyclic.
\end{corollary}




\section*{acknowledgments}

The authors thank E. Pardo and M. Siles Molina  for their valuable discussions
during the preparation of this paper.  The authors also thank the referee for a very careful reading of, and suggested changes to, the initial version of the manuscript.



\begin{thebibliography}{}



\bibitem{AALP} \textsc{Abrams, G.,  \'{A}nh, P.N.,  Louly, A.,  Pardo, E.}:  The classification question for Leavitt path algebras, \emph{J. Algebra}, to appear.  ArXiV: 0706.3874.

\bibitem{AA3}  \textsc{Abrams, G.,  Aranda Pino, G.}: The Leavitt path algebra of arbitrary graphs, \emph{Houston J. Math}, to appear.


\bibitem{AAS} \textsc{Abrams, G., Aranda Pino, G., Siles Molina, M.}: Finite-dimensional
Leavitt path algebras, \emph{J. Pure Appl. Algebra.} \textbf{209(3)},  753--762 (2007).

\bibitem{AM} \textsc{\'{A}nh, P.N., M\'{a}rki, L.}: Morita equivalence for rings without identity,
\emph{Tsukuba J. Math.} \textbf{11(1)}, 1-16 (1987).

\bibitem{AMP} \textsc{Ara, P., Moreno, M.A., Pardo, E.}:, Nonstable K-Theory
for graph algebras, \emph{Algebra Represent. Theory} \textbf{10(2)}, 157--178  (2007).


\bibitem{CK}  \textsc{Camillo, V., Khurana, D.}: A characterization of unit regular rings, \emph{Comm. Algebra} \textbf{29(5)}, 2293--2295   (2001).

\bibitem{CY}  \textsc{Camillo, V., Yu, H.-P.}: Stable range one for rings with many idempotents, \emph{Trans. A.M.S.} \textbf{347(8)}, 3141--3147  (1995).


\bibitem{G}  \textsc{Goodearl, K.}:  Leavitt path algebras and direct limits, to appear.   ArXiV:0712.2554v1

\bibitem{G2}  \textsc{Goodearl, K.}: {\it Von Neumann Regular Rings}.  Krieger Publ., Malabar, FL (1991).  ISBN 0-89464-632-X.

\bibitem{K}  \textsc{Kaplansky, I}: Topological representation of algebras.  II, \emph{Trans. A.M.S.} \textbf{68(1)}, 62--75  (1950).

\bibitem{NZ}  \textsc{Nicholson, W.K, Zhou, Y.}: Clean Rings: a survey.    In: {\it Advances in Ring Theory: Proceedings of the 4th China-Japan-Korea International Conference}, pp. 181--198.  World Sci. Publ., Hackensack, N.J.  (2005).  ISBN: 981-256-425-X.

\bibitem{R}  \textsc{Rowen, L.}: Examples of semiperfect rings, \emph{Israel J. Math.} \textbf{65(3)}, 273--283  (1989).

\bibitem{RS}  \textsc{Raeburn, I., Szyma\'{n}ski, W.}: Cuntz-Krieger algebras of infinite graphs and matrices, \emph{Trans. A.M.S.} \textbf{356(1) }, 39--59  (2003).





\end{thebibliography}


\end{document}